\documentclass[12pt]{amsart}
\begin{document}
\newcommand{\lar}{\longrightarrow}
\newcommand{\q}{^{[q]}}

\newtheorem{theorem}{Theorem}[section]

\newtheorem{prop}[theorem]{Proposition}
\newtheorem{lemma}[theorem]{Lemma}
\newtheorem{corollar}[theorem]{Corollary}
\newtheorem*{Corollary}{Corollary}
\newtheorem{obs}[theorem]{Observation}
\theoremstyle{definition}

\newtheorem{note}[theorem]{Note}
\newtheorem{definition}[theorem]{Definition}
\newtheorem{example}[theorem]{Example}
\newtheorem*{notation}{Notation}
\newtheorem*{conj}{Conjecture}
\newtheorem*{claim}{Claim}
\newtheorem*{question}{Question}

\newcommand{\li}{\tilde}
\newcommand{\aaa}{\mathfrak{a}}
\newcommand{\bbb}{\mathfrak{b}}
\newcommand{\ccc}{\mathfrak{c}}
\newcommand{\ub}{\underline{b}}

\newcommand{\m}{\mathfrak{m}}

\newcommand{\bs}{\boldsymbol}
\newcommand{\ld}{\ldots}
\newcommand{\cd}{\cdots}

\newcommand{\ui}{\underline{i}}
\newcommand{\und}{\underline}
\newcommand{\arr}{\Rightarrow}
\newcommand{\inc}{\subseteq}
\newcommand{\w}{\omega}
\newcommand{\stI}{{}^\aaa I^*}
\newcommand{\st}{{}^\aaa}
\newcommand{\Tor}{\mathrm{Tor}^R}
\newcommand{\Ext}{\mathrm{Ext}_R}
\newcommand{\x}{\underline{x}}
\newcommand{\ol}{\overline}
\newcommand{\syz}{\mathrm{syz}}
\newcommand{\Soc}{\mathrm{Soc}}

\title{A new version of $\mathfrak{a}$-tight closure}
\author{Adela Vraciu}
\address{Department of Mathematics, University of South Carolina, Columbia SC 29205}
\email{vraciu@math.sc.edu}
\subjclass{13A35}
\thanks{The author thanks N. Epstein, A. Kustin, and K. Yoshida for useful conversations on the topic of this manuscript. The author was partially supported by a grant from the University of South Carolina Research Foundation.}
\date{March 30, 2007}


\begin{abstract}
Hara and Yoshida introduced a notion of $\aaa$-tight closure in 2003, and they proved that the test ideals given by this operation correspond to multiplier ideals. However, their operation is not a true closure. The alternative operation introduced here is a true closure. Moreover, we define a joint Hilbert-Kunz multiplicity that can be used to test for membership in this closure. We study the connections between the Hara-Yoshida operation and the one introduced here, primarily from the point of view of test ideals. We also consider variants with positive real exponents.
\end{abstract}
\maketitle

\section{Introduction}

In \cite{HY}, Hara and Yoshida introduced a notion of $\aaa$-tight closure, which generalizes the ``classical'' tight closure of Hochster and Huneke introduced in \cite{HH1}. The main motivation for their work is the connection between the test ideals given by this operation and multiplier ideals, which generalizes previous results of Hara (\cite{Ha}) and Smith (\cite{Sm}). The advantage of this version of test ideal is that it allows them to study multiplier ideals for pairs, not only the multiplier ideal of a variety.

However, the Hara-Yoshida $\aaa$-tight closure is not a true closure operation, since it gets (potentially) larger when iterated. The version introduced in this paper is a true closure, and it is always contained in the Hara-Yoshida $\aaa$-tight closure. We establish several other connections between these operations. We prove that for a Gorenstein graded algebra of dimension at least 2, the test ideals given by these two operations are the same (Theorem \ref{testideals}). The Hara-Yoshida $\aaa$-tight closure of an ideal $I$ is denoted $I^{*\aaa}$, while the new version introduced here will be denoted ${}^{\aaa}I^*$.

We  define a joint Hilbert-Kunz multiplicity associated to two $\m$-primary ideals $\aaa$ and $I$, and we prove that this multiplicity can be used to test for membership in our version of $\aaa$-tight closure.
This is similar to the way in which the Hilbert-Samuel multiplicity is used to test for membership in the integral closure, and the Hilbert-Kunz multiplicity is used to test for membership in tight closure. 

There are versions of both closures, as well as of the joint multiplicity, in which positive real numbers are allowed as exponents. For fixed ideals $I$ and $\aaa$, we study the question of how $I^{*\aaa^t}$ and ${}^{\aaa^t}I^*$ vary with $t$.
This question is related to the notion of jumping exponents (in characteristic zero), or F-thresholds (in positive characteristic).

In this paper, $R$ will denote a Noetherian ring of positive characteristic $p>0$ and Krull dimension $d>0$, and $q=p^e$ will always denote a power of the characteristic. Most of the time, $R$ will be assumed to be either local or graded. $R^o$ is the set of elements in $R$ that are not in any minimal prime  of $R$. If $I\subset R$ is an ideal, $I\q$ denotes the ideal generated by all $i^q$, when $i \in I$.

\section{Definitions and elementary properties}

\begin{definition}(\cite{HY})
Let $\aaa, I$ be ideals in $R$, and $x \in R$. 
We say that $x \in I^{*\aaa}$ if there exists $c \in R^0$ such that $c\aaa ^q x^q \subseteq  I\q$ for all $q=p^e\gg 0$.
\end{definition}
 
\begin{definition}
Let $\aaa, I$ be ideals in $R$, and $x \in R$. 
We say that $x \in {}^{\aaa}I^*$ if there exists $c \in R^0$ such that $c\aaa ^q x^q \subseteq \aaa ^q I\q$ for all $q=p^e\gg 0$.
\end{definition}

In the case when $\aaa=R$, both of the above definitions coincide with the definition of the usual tight closure of Hochster and Huneke (\cite{HH1}), which is denoted $I^*$.
Some elementary properties of these operations are summarized below.
\begin{obs}\begin{enumerate}

\item For all $\aaa$ and $I$, we have $I^* \subseteq {}^{\aaa} I^* \subseteq I^{*\aaa}$.

\item If $\aaa =(f)$ is a principal ideal, then $I^{*\aaa}=I^* : f$. In particular, $(I^{*(f)})^{*(f)} \ne I^{*(f)}$ when $(R, \m)$ is local, $I$ is $\m$-primary, and $f \in \m \, \backslash \, I$.

\item If $\aaa =(f)$ is a principal ideal, and $f$ is a non-zerodivisor on $R$, then
${}^{\aaa}I^*=I^*$.

\item For all $\aaa$ and $I$, $\st (\stI)^*=\stI$. In other words, $\stI$ is a true closure operation.
\end{enumerate}
\end{obs}
\begin{proof}
(1), (2) and (3) are trivial.

(4) Let $x \in \st (\stI)^*$. Then there exists $c \in R^o$ such that
$c\aaa^q x^q \subseteq \aaa^q (\stI)\q$ for all $q=p^e$.
Also, there exists $c'\in R^o$ such that 
$c'\aaa^q (\stI)\q \subseteq \aaa^q I\q$. 
Combining these two inclusions, we get
$cc'\aaa^q x^q\subseteq \aaa^q I\q$.
\end{proof}

The notion of {\it test element for tight closure} was defined in \cite{HH1}, and it was proved that test elements (for tight closure) exist in very general classes of rings.
\begin{definition}
An element $c \in R^o$ is called a {\it test element} for tight closure if we have $cx \in I$ for every ideal $I$ and every $x \in I^*$.
\end{definition}
  \begin{prop}\label{briancon}
Assume that $R$ has test elements for the usual tight closure. Then for any ideals $\aaa$ and $I$, with $I$ of positive height, we have
${}^{\aaa}I^* \subseteq \overline{I}$.
\end{prop}
\begin{proof}
 By the usual determinant trick, $cx^q \aaa^q \subseteq \aaa^q I\q$ implies $cx^q \in \overline{I\q} \subseteq \overline{I^q}\subseteq (I^{q-n+1})^*$, where $n$ is the minimal number of generators of $I$. The last inclusion is by the tight closure version of the Brian\c{c}on-Skoda theorem (\cite{HH1}, Theorem 5.4). Let $d\in R^o$ be a test element, and $f \in I^{n-1}\cap R^o$ a fixed element. Then we have
$cdf x^q \in I^q$ for all $q =p^e$, which shows that $x \in \overline{I}$. 
\end{proof}

Versions in which positive real numbers occur as exponents can be defined for both operations:
\begin{definition}
Let $\aaa, I\subset R$ be ideals, and let $t>0$ be a real number. Let $x \in R$. For any real number $r$, $\lceil r\rceil $ denotes the smallest integer greater than or equal to $r$.
\begin{enumerate}

\item We say that $x \in I^{*\aaa^t}$ if there exists $c \in R^o$ such that $cx^q \aaa^{\lceil tq \rceil } \subseteq I\q$ for all $q=p^e$.

\item We say that $x \in {}^{\aaa ^t}I^*$ if there exists $c \in R^o$ such that $cx^q \aaa^{\lceil tq \rceil }\subseteq \aaa^{\lceil tq \rceil }I\q$ for all $q=p^e$.
\end{enumerate}
 \end{definition}

Each of these operations gives rise to a test ideal as follows:
\begin{definition}
Let $\aaa\subseteq R$ be an ideal, and let $t>0$ be a real number.
We define
$$
\tau (\aaa^t):=\bigcap (I: I^{*\aaa^t})\ \ \ \ T_{\aaa^t}:= \bigcap (I:{}^{\aaa^t}I^*),
$$
where each intersection ranges over all the ideals $I\subseteq R$.
\end{definition}

Part (1) of the next observation was noted in \cite{HT}, where it was used to prove Skoda's theorem for test ideals. Part (2) is an analogue for the new $\aaa$ -tight closure. 
\begin{obs}
Let $\aaa, I \subset R$ be ideals, and $t > \nu (\aaa)$ a real number, where $\nu(\aaa)$ denotes the minimal number of generators of $\aaa$ . Then:
\begin{enumerate}

\item $I^{*\aaa^t}=(I^{*\aaa^{t-1}}): \aaa$.

\item ${}^{\aaa^t}I^*=({}^{\aaa^{t-1}}(\aaa I)^*) : \aaa$. 
\end{enumerate}
\end{obs}
\begin{proof}
(1) The proof of this statement can be found as part of the proof of Theorem 4.1 in \cite{HT}.

(2) Note that we have $\aaa ^r = \aaa ^{[q]} \aaa^{r-q}$ for all $r> \nu (\aaa) q$. 
We have
$x \in {}^{\aaa^t}I^* \Leftrightarrow cx^q \aaa^{\lceil tq\rceil } \subseteq \aaa^{\lceil tq\rceil }I\q \Leftrightarrow 
c (\aaa x)\q \aaa^{\lceil tq \rceil -q }\subseteq \aaa^{\lceil tq \rceil -q} (\aaa I ) \q$ and the conclusion follows since $\lceil tq \rceil -q = \lceil (t-1) q\rceil $.
\end{proof}
We establish two connections between the two versions of $\aaa$-tight closure. The first result, Proposition~\ref{highdeg}, shows that for elements of large enough degree in a graded ring, membership in one of these closures is equivalent to membership in the other. The second result, Proposition~\ref{intermediary} shows that, under certain assumptions, every element in the Hara-Yoshida $\aaa$-tight closure must satisfy a stronger condition, which bridges the gap between the Hara-Yoshida definition and the one introduced in this paper.

We establish the following notation, which will be in effect throughout this paper when graded rings are considered.
\begin{notation}\label{Notation}
If $R$ is a finitely generated graded algebra over a field,
$R=\oplus _{n \ge 0}R_n$, we will denote $R_{_+} := \oplus _{n >0} R_n$ the unique maximal homogeneous ideal of $R$. We will let $y_1,  \ldots, y_s$ be algebra generators for $R$, and let $\beta_1, \ldots, \beta_s$ be their degrees. Set $\beta = \mathrm{max}(\beta_i), \beta'=\mathrm{min}(\beta_i)$.

We say that $R$ is standard graded if $\beta _i =1$ for all $i$.
\end{notation}
\begin{prop}\label{highdeg}
Let $R$ be a finitely generated graded algebra over a field  and let $\aaa \subset R$ be a homogeneous $R_{_+}$-primary ideal, so that $R_{_+}^k \subseteq \aaa \subseteq R_{_+}^l$ for some integers $l \le k$. Let $I=(f_1, \ldots, f_n)$ be a homogeneous ideal, and $x \in R_N$ with $N \ge \beta k-\beta 'l + \mathrm{max}(\mathrm{deg}(f_i))$. Then $x \in I^{*\aaa} \Leftrightarrow x \in {}^{\aaa}I^*$.
\end{prop}
\begin{proof}
Assume that $x \in I^{*\aaa}$ and $\mathrm{deg}(x) \ge \beta k-\beta 'l + \mathrm{max}(\mathrm{deg}(f_i))$. For each homogeneous $h\in \aaa ^q$, we have $\mathrm{deg}(h)\ge \beta 'lq$. We can write
$cx^q h=\Sigma _{i=1}^n a_i f_i^q$ with $c \in R^o$, $a_i \in R$ homogeneous elements, so that
$\mathrm{deg}(a_i)=\mathrm{deg}(c)+ q\mathrm{deg}(x)+ \mathrm{deg}(h)-q\mathrm{deg}(f_i)\ge \beta kq$ for each $i$. Thus, $a_i \in R_{\ge \beta kq}\subseteq R_+^{kq}\subseteq \aaa^q$ (the first inclusion follows because any element in $R_{\ge \beta kq}$ can be written as a linear combination of monomials $y_1^{i_1} \cdots y_s^{i_s}$ with $i_1\beta _1 + \ldots + i_s \beta _s \ge \beta kq$, which implies that $i_1 + \ldots + i_s\ge kq$).
\end{proof}

\begin{obs}\label{isic}
If $R$ is standard graded, so $\beta =\beta'$, $\aaa=R_{_+}^r$, and all the generators of $I$ have the same degree, then we have ${}^{\aaa}I^*=I^{*\aaa} \cap R_{\ge N}$, where $N$ denotes the common degree of the generators of $I$. 
\end{obs}
\begin{proof}
Let $x \in {}^{\aaa}I^*$, so there exists $c \in R^o$ (which can be assumed homogeneous) such that $cx^q \aaa^q \subseteq \aaa^q I\q$. Taking degrees of both sides yields $\mathrm{deg}(c) + q \mathrm{deg}(x) + qr \ge qr + qN$, so that $\mathrm{deg}(x) \ge N$. This shows that  ${}^{\aaa}I^*\subseteq I^{*\aaa} \cap R_{\ge N}$. The other inclusion is contained in Proposition~\ref{highdeg}.
\end{proof}
This observation might suggest that ${}^{\aaa}I^*=I^{*\aaa} \cap \overline{I}$ for $\m$-primary ideals $I$. This is in fact not true (however, ${}^{\aaa}I^*\subseteq  I^{*\aaa}\cap \overline{I}$ is always true), as seen in the following example.
\begin{example}\label{notic}
Let $R=k[x, y]$, $I=(x^2, y^4)$, $\aaa=(x, y)^3$. Then we have $I^{*\aaa}=I+(xy^2, y^3)=I:(x, y)^2$,
$I^{*\aaa} \cap \overline{I}=I+(xy^2)$, and
${}^{\aaa}I^*=I+(xy^3)=I:(x, y)$.
\end{example}
\begin{proof}
If $i + j \ge 3q$, we have $i \ge q$ or $j \ge 2q$. In either case we have
$x^iy^j x^{q}y^{2q}\in(x^{2q}, y^{4q})$, and thus $xy^2 \in I^{*\aaa}$. Similarly, if $i+j\ge 3q$ we have $i\ge 2q$ or $j \ge q$; in either case, $x^iy^jy^{3q}\in (x^{2q}, y^{4q})$, and thus $y^3\in I^{*\aaa}$. Also note $(xy^2)^2 \in I^2$, so $xy^2 \in \overline{I}$. However, $y^3\notin \overline{I}$ (one can see this from the Newton polygon, for instance).

To see that $xy^2 \notin {}^{\aaa}I^*$, we prove the stronger fact that $xy^2 \notin{}^{(x, y)^n}I^*$ for any $n \ge 3$. This will suffice to prove the last statement, since all the ideals under consideration are monomial. Assume the contrary, so that there exists $c \in R^o$ such that
$cx^iy^j x^qy^{2q}\in (x, y)^{nq}I\q$ for all $i, j$ with $i+j=nq$, for some $n$.

Choose $i=\displaystyle \lceil \frac{q}{2} \rceil , j =(n-1) q + \lfloor \frac{q}{2} \rfloor$.
We obtain $cx^{\lceil 3q/2 \rceil } y^{(n-3) q + \lfloor q/2\rfloor } y^{4q}=ax^{2q}+b y^{4q}$ with $a, b \in (x, y)^{nq}$. This is clearly impossible since the degree of $x^{\lceil 3q/2 \rceil } y^{(n-3) q+ \lfloor q/2 \rfloor }$ is $(n-1)q$, and the degree of $c$ is a constant.
\end{proof}

\begin{prop}\label{intermediary}
Let $(R, \m)$ be an excellent normal domain such that its completion is a domain. Let $I, \aaa \subset R$ be ideals, and assume that $\aaa $ is not a principal ideal. Then there exists a $Q_0=p^{e_0}$ and a $c\in R^o$ such that for all $x \in I^{*\aaa}$ we have $cx^q \aaa ^q \subseteq \m^{q/Q_0}I\q$ for all $q\gg 0$.
\end{prop}
Note that if $\aaa $ is $\m$-primary, then we can replace $\m^{q/Q_0}$ by $\aaa ^{q/Q_0}$ by choosing a possibly larger $Q_0$.
\begin{proof}
First note that there is no loss of generality in assuming that $I$ is $*$-independent, i.e. $I=(f_1, \ldots, f_n)$ with $f_i \notin (f_1, \ldots, \hat{f_i}, \ldots, f_n)^*$ for all $i$. That is because one can find a $*$-independent $I_0 \subseteq I$ with $I_0^* = I^*$ (by omitting generators of $I$ that are redundant up to tight closure), and it is easy to see that $I_0^* = I^*$ implies $I_0^{*\aaa}=I^{*\aaa}$.

Let $\aaa =(a_1, \ldots, a_s)$, with $s \ge 2$, and $I=(f_1, \ldots, f_n)$. The $*$-independence assumption implies that there exists $q_1$ such that 
$$(f_1^q, \ldots, \hat{f_i}^q, \ldots, f_n^q):f_i ^q \subseteq \m ^{[q/q_1]}$$ for all $q$ and all $i$ (cf. Proposition 2.4 in \cite{Ab}).

Since $R$ is normal, we have $a_l \notin \overline{(a_k)}$ for any $1\le k \ne l\le s$, and we can choose $q_2 \gg 0$ such that $a_l \notin \overline{(a_k, \m^{q_2/q_1})}$. Also choose $q_2 \ge s$. In particular, $a_l \notin (a_k, \m^{q_2/q_1})^*$ and thus we can choose $q_0$ such that $(a_k^q, \m^{[qq_2/q_1]}): a_l^q \subseteq \m^{q/q_0}$ (using Proposition 2.4 in \cite{Ab} again).

 We have $cx^{qq_2}\aaa^{qq_2} \subseteq I^{[qq_2]}$ for a fixed $c \in R^o$ and all $q$.
 Fix an element $a_1^{i_1} \cdots a_s^{i_s} \in \aaa^{qq_2}$ and write
$cx^{qq_2} a_1^{i_1} \cdots a_s^{i_s}= b_1f_1^{qq_2} + \ldots + b_n f_n^{qq_2}$. The choice of $q_2$ guarantees that $i_k \ge q$ for some $k$. Choose an index $l \ne k$ and consider the element $a_1^{j_1}  \cdots a_s^{j_s} \in \aaa^{qq_2}$ with $j_k = i_k -q$, $j_l = i_l +q$, and $j_{\tau } = i_{\tau }$ for all other $\tau = 1, \ldots, s$.
We have
$cx^{qq_2}a_1^{j_1}  \cdots a_s^{j_s}=b_1'f_1^{qq_2} + \ldots b_n'f_n^{qq_2}$.
Multiplying the first equation by $a_l^q$ and the second equation by $a_k^q$ yields $(b_i a_l^q - b_i' a_k^q) \in (f_1^{qq_2}, \ldots, \hat{f_i^{qq_2}}, \ldots, f_n^{qq_2}) :f_i^{qq_2}\subseteq \m^{[qq_2/q_1]}$, and therefore 
$b_i \in (a_k^q, \m^{[qq_2/q_1]}): a_l^q\subseteq \m^{q/q_0}$. 
This holds for all $i =1, \ldots, n$, and for any choice of the multi-index $(i_1, \ldots, i_s)$. We get the desired conclusion by choosing $Q_0 =q_2q_0$.
\end{proof}

\section{Joint Hilbert-Kunz multiplicities}
The idea of associating a multiplicity to a pair or more ideals (the so-called mixed multiplicity) first appeared in \cite{Bt}, and the notion was extensively studied by many other authors, including B. Tessier, D. Rees and I. Swanson. The idea of a multiplicity coming from length functions involving both ordinary and Frobenius powers can be found in work of Hanes (\cite{Hn}). The joint Hilbert-Kunz multiplicity introduced here bares a resemblance to each of these previous multiplicities, but is different from them.

Assume that $(R, \m)$ is local, let $I, \aaa \subset R$ be $\m$-primary ideals, $M$ a finitely generated $R$-module, and $t>0$ a fixed real number.

We study the function
$$
\ell_M(q):=\lambda \left( \frac{M}{\aaa ^{\lceil qt \rceil} I\q M}\right),
$$
where $q=p^e$. Note that $\aaa ^{\lceil qt \rceil}$ is an ordinary power where the exponent is obtained by taking the least integer which is greater than or equal to $tq$, while $I\q$ is a Frobenius power. We will write $\ell(q)$ for $\ell _R(q)$.

\begin{theorem}\label{hk}
Let $R, I, \aaa, M, t$ be as above, and let $d$ be the Krull dimension of $R$. Then there is a $c>0$ such that
$$
\ell_M(q)=cq^d+\mathcal{O}(q^{d-1}).
$$
\end{theorem}
If $M=R$, we call $c$ the {\it mixed Hilbert-Kunz multiplicity} of the pair $(\aaa^t, I)$ and we denote it $e_{HK}(\aaa^t, I)$.

The proof of the Theorem follows essentially the same steps as in Monsky's paper (\cite{Mo}). We will follow closely the outline of his paper.
\begin{lemma}\label{L1}
Assume that there is an $h\in R^o$ such that $hM=0$. Then there exists $a>0$ such that $\ell_M(q)\le aq^{d-1}$.
\end{lemma}
\begin{proof}
Let $n$ be the number of generators of $I$. Then we have $I\q \supseteq I^{nq}$. Also, $\lceil t \rceil q \ge \lceil tq \rceil$, so $\aaa ^{\lceil tq \rceil } \supseteq  \aaa^{\lceil t \rceil q}$, and it follows that
$$
\ell_M(q)\le \lambda\left( \frac{M}{(\aaa^{\lceil t \rceil } I^n)^q M}\right),
$$ which is a Hilbert-Samuel function over the ring $R/h$, a ring of Krull dimension at most $d-1$, and thus it is bounded by $aq^{d-1}$ for some $a>0$.
\end{proof}

\begin{lemma}\label{L2}
Let $M, N$ be finitely generated $R$-modules such that $M_{p_i}\cong N_{p_i}$ for every minimal prime $p_i$ of $R$. Then
$|\ell_M(q)-\ell_N(q)|\le \mathcal{O}(q^{d-1})$.
\end{lemma}
\begin{proof}
Let $S=R\, \backslash \, \bigcup p_i$. We have $S^{-1}M\cong S^{-1}N$. Since $S^{-1}\mathrm{Hom}_R(M, N)\cong \mathrm{Hom}_{S^{-1}R}(S^{-1}M, S^{-1}N)$, we have a homomorphism 
$\phi : M\rightarrow N$ such that $S^{-1}\phi $ is bijective. Unlocalizing, we get an element $h \in S$ such that $h$ annihilates the cokernel $C$ of $\phi$.
Consider the exact sequence
$$
\frac{M}{\aaa^{\lceil tq\rceil}  I\q M}\lar \frac{N}{\aaa^{\lceil tq \rceil } I\q N}\lar \frac{C}{\aaa^{\lceil tq \rceil} I\q C}\lar 0$$
Lemma ~\ref{L1} gives 
$$
\ell_N(q)-\ell_M(q)\le \ell_C(q)\le  aq^{d-1}$$
for some $a>0$.
Now repeat the argument with the roles of $M, N$ reversed in order to get
$$
\ell_M(q)-\ell_N(q)\le bq^{d-1}$$
for some $b>0$.
\end{proof}

\begin{definition}
Let $M_{(e)}$ be $M$ viewed as an $R$-module via the Frobenius map $F^e: R\lar R$.
Note that $\underline{ }_{(e)}$ is an exact functor, and, if we assume that the residue field of $R$ is perfect, we have  
$$
\ell_{M_{(e)}}(q)= \lambda\left(\frac{M}{(\aaa ^{\lceil tq\rceil })^{[p^e]}I^{[qp^e]}M}\right).
$$
\end{definition}

The following is the one essential ingredient we need in addition to Monsky's ideas:

\begin{obs}\label{computation}
Let $R, \aaa, I, t$ be as above, $e>0$ a fixed integer. By prime avoidance, we can choose generators $f_1, \ldots, f_n$ of $\aaa$ that are in $R^o$. Let $f:=f_1 \cdots f_n$.
Then:

{\bf a.}  $\aaa^{\lceil tqp^e\rceil }\subseteq (\aaa^{\lceil tq\rceil})^{[p^e]} : f^{p^e}$.

{\bf b.} Assume that the residue field of $R$ is perfect. Then
$| \ell_{M_{(e)}}(q) -  \ell _M (p^e q )| \le \mathcal{O}(q^{d-1}).$
\end{obs}
\begin{proof}
{\bf a.} The generators of $\aaa ^{\lceil tqp^e\rceil }$ are of the form
$F=f_1^{a_1p^e+i_1}\cdots f_n^{a_np^e+i_n}$, where $0\le i_k\le p^e-1$ for all $k$, and $$(a_1+\ldots +a_n)p^e+i_1+\ldots +i_n=\lceil tqp^e\rceil \ge (\lceil tq \rceil-1)p^e ,$$ with all $a_k, i_k \in {\bf Z}$. It follows that
$a_1+\ldots +a_n \ge \lceil tq\rceil -1-n+n/p^e$. Since $a_k \in {\bf Z}$ for all $k$, it must be that $a_1+\ldots +a_n \ge \lceil tq\rceil -n$ and thus $(a_1+1)+\ldots +(a_n+1)\ge \lceil tq\rceil $, and so $f_1^{p^e-i_1}\cdots f_n^{p^e-i_n}F=(f_1^{a_1+1}\cdots f_n^{a_n+1})^{p^e} \in (\aaa ^{\lceil tq\rceil} )^{[p^e]}$. 

{\bf b.} We have
$$
\ell_{M_{(e)}}(q) -\ell_M(p^eq)=
\lambda \left( \frac{\aaa^{\lceil tqp^e \rceil} I^{[qp^e]} M}{
(\aaa^{\lceil tq \rceil })^{[p^e]}I^{[qp^e]}M}\right) \le
$$
$$\lambda\left( \frac{(\aaa^{\lceil tq \rceil })^{[p^e]}I^{[qp^e]}M: f^{p^e}}{(\aaa^{\lceil tq \rceil })^{[p^e]}I^{[qp^e]}M}\right)
=
\lambda \left(\frac{M}{(\aaa^{\lceil tq \rceil })^{[p^e]}I^{[qp^e]}M+(f^{p^e})M}\right).
$$

The inequality above follows from part {\bf a}. The second equality follows from the general fact that
for any $\m$-primary ideal $J\subset R$, and any element $g\in R$, we have
$$
\lambda \left( \frac{JM: g}{JM}\right) =\lambda \left( \frac{M}{(J, g)M}\right)$$
applied to $J=(\aaa^{\lceil tq \rceil })^{[p^e]}I^{[qp^e]}$ and $g=f^{p^e}$.
(Proof of the general fact: consider the short exact sequence
$$0 \lar \frac{M}{JM:g}\lar \frac{M}{JM} \lar \frac{M}{JM+(g)M}\lar 0$$
where the first map is multiplication by $g$.) 
Lemma ~\ref{L1} now gives the desired conclusion, since 
$$\lambda \left(\frac{M}{(\aaa^{\lceil tq \rceil })^{[p^e]}I^{[qp^e]}M+(f^{p^e})M}\right)$$
 is a joint Hilbert-Kunz function over the $d-1$ dimensional ring $R/(f^{p^e})$.
\end{proof}

\begin{lemma}\label{L3} Assume that the residue field of $R$ is perfect.
Let $$0\lar M'\lar M\lar M''\lar 0$$ be a short exact sequence of finitely generated $R$-modules.
Then we have
$$
\ell_M(q)=\ell_{M'}(q)+\ell_{M''}(q)+\mathcal{O}(q^{d-1}).
$$
\end{lemma}
\begin{proof}
{\bf Case 1:} Assume that $R$ is reduced. For each minimal prime ${p_i}$ of $R$, $R_{p_i}$ is a field and it follows that
$M_{p_i}\cong (M'\oplus M'')_{p_i}$. The conclusion follows from Lemma ~\ref{L2}.

{\bf Case 2:} Let $\mathfrak{n}$ denote the nilradical of $R$, and choose $e$ such that $\mathfrak{n}^{[p^e]}=0$. Note that $M_{(e)}$ is annihilated by $\frak{n}$ for every module $M$.
We get a short exact sequence of $R/\frak{n}$ modules:
$$
0\lar M'_{(e)} \lar M_{(e)} \lar M''_{(e)} \lar 0,$$
and now we can apply the result from case 1 in conjunction with Obs. ~\ref{computation}.
\end{proof}

\begin{lemma}\label{L4}
Assume that $R$ is a domain with perfect residue field. Then there exists $c>0$ such that
$$
\ell(q)=cq^d+\mathcal{O}(q^{d-1}).
$$
\end{lemma}
\begin{proof}
It is known that the rank of $R_{(1)}$ as an $R$-module is $p^d$.

Apply Lemma \ref{L2} to the $R$-modules $R_{(1)}$ and $R^{p^d}$. We get
$$
|\lambda \left(\frac{R}{(\aaa^{\lceil tq\rceil} )^{[p]}I^{[pq]}}\right) -p^d\ell(q)|=|\ell_{R_{(1)}}(q)-p^d\ell(q)|\le a'q^{d-1}$$
for some $a'>0$, and by Obs. ~\ref{computation} it follows that 
$$
|\ell(pq)-p^d \ell_R(q) | \le a q^{d-1}
$$ for some $a$.
Thus, we have 
$$
\left | \frac{\ell(pq)}{(pq)^d} - \frac{\ell(q)}{q^d}\right| \le \frac{a'}{p^dq}.
$$
It follows that
$$
\left | \frac{\ell(q'q)}{(q'q)^d} - \frac{\ell(q)}{q^d}\right| \le \frac{a'}{p^dq} \frac{1-\frac{1}{q'}}{1-\frac{1}{p}},
$$
thus $\{\ell(q)/q^d\}$ is a Cauchy sequence. Let $c:=\lim_{q\rightarrow \infty} \ell(q)/q^d$.
If we keep $q$ fixed and let $q'\rightarrow \infty$, we get 
$$
|\ell(q)/q^d-c|\le \frac{\alpha '}{q}
$$
for some $\alpha '$ and all $q$, and thus
$
|\ell(q)-cq^d|\le \alpha 'q^{d-1},
$ 
or in other words $\ell(q)=cq^d+\mathcal{O}(q^{d-1})$.
\end{proof}

Now we are ready to prove the general case of Theorem~\ref{hk}.
\begin{proof}
Since every finitely generated module $M$ has a filtration $(0)=M_0 \subset M_1 \subset \ldots \subset M_n =M$ with quotients $M_{i+1}/M_i \cong R/P_i$, with $P_i$ prime ideals, the general case follows from Lemma ~\ref{L4} by repeated application of lemma ~\ref{L3}.\

In order to remove the assumption that the residue field is perfect, note that length is preserved by faithfully flat base change. Thus, we can pass to completion, so that
$R$ is a quotient of a formal power series ring $K[[X_1, \ldots, X_n]]$, and we can replace $R$ by $R\otimes _K F$, where $F$ is an algebraic closure of $K$.
\end{proof}

\begin{lemma}\label{realmult}
Let $\aaa \subset R$ be an $\m$-primary ideal and $t>0$ a real number.
Then
$$
\lim_{q\rightarrow \infty} \lambda (\frac{R}{\aaa^{\lceil tq \rceil }})/q^d=\frac{t^d e(\aaa)}{d!}
$$
\end{lemma}
\begin{proof}
First note that there exists a sequence of rational numbers $\{ k_n /q_n \}$ with denominators of the form $q_n =p^{e_n}$ such that
$\displaystyle \frac{k_n}{q_n} \le  t < \frac{k_n+1}{q_n}$, and $q_n < q_{n+1}$, so that $t =\lim_{n \rightarrow \infty } q_n /k_n$.
For instance, take $q_n=p^n$, $k_n=\lfloor tp^n\rfloor$.

For $n$ fixed and $q=p^e\gg 0$, we have $\displaystyle k_n\frac{q}{q_n} \le \lceil tq \rceil \le (k_n+1) \frac{q}{q_n}$, and 
$$
\lambda \left( \frac{R}{\aaa^{(k_n +1)q/q_n}}\right)/q^d=
e(\aaa) \frac{(k_n+1)^d (q/q_n)^d}{d!} + \mathcal{O}(q^{d-1}),
\ \mathrm{and}\, 
$$
$$\lambda \left( \frac{R}{\aaa^{k_nq/q_n}}\right)/q^d=e(\aaa)\frac{(k_n)^d (q/q_n)^d}{d!} + \mathcal{O}(q^{d-1}),
$$
so for all $n$ we have
$$
e(\aaa) \frac{(k_n)^d}{d! q_n)^d} \le \lim_{q\rightarrow \infty} \lambda\left( \frac{R}{\aaa^{\lceil tq \rceil } }\right)/q^d \le e(\aaa)\frac{(k_n+1)^d}{d!q_n)^d}
 $$
and the desired result follows by taking the limit when $n \rightarrow \infty$.
\end{proof}
\begin{theorem}\label{continuity}
If $\aaa, I$ are fixed $\m$-primary ideals, then $e_{HK}(\aaa^t, I)$ is a continuous function of $t$.
\end{theorem}
\begin{proof}
Let $t< t'$ be positive real numbers. Then
$$e_{HK}(\aaa^{t'}, I) - e_{HK}(\aaa^t, I)=
\lim_{q\rightarrow \infty }\frac{ \lambda \left(\aaa^{\lceil tq\rceil }I\q/\aaa^{\lceil t'q \rceil }I\q \right)}{q^d}.
$$
Let $I=(f_1, \ldots, f_n)$. Then we have a composition series
$$
\aaa^{\lceil t'q \rceil }I\q \subseteq \aaa^{\lceil t'q \rceil }I\q + f_1^q \aaa^{\lceil tq \rceil} \subseteq \ldots \subseteq \aaa^{\lceil t'q \rceil }I\q+ (f_1^q, \ldots, f_i^q) \aaa^{\lceil tq \rceil}\subseteq  \ldots \subseteq  \aaa^{\lceil tq \rceil}I\q
$$
For $i=1, \ldots, n$, 
let $K_i = \aaa^{\lceil t'q \rceil }I\q + (f_1^q, \ldots, f_{i-1}^q) \aaa^{\lceil tq \rceil}$.
Then we have
$$
\lambda \left(\frac{\aaa^{\lceil tq \rceil}I\q}{\aaa^{\lceil t'q \rceil }I\q}\right)=
\Sigma_{i=1}^n \lambda \left( \frac{K_i+ \aaa^{\lceil tq \rceil}f_i^q }{K_i}\right)
$$
Note that
 we have a surjective map given by multiplication by $f_i^q$:
$$
\frac{\aaa^{\lceil tq \rceil}}{(K_i : f_i^q)\cap \aaa^{\lceil tq \rceil}} \rightarrow \frac{\aaa^{\lceil tq \rceil}f_i^q}{K_i \cap \aaa^{\lceil tq \rceil}f_i^q}\cong \frac{K_i + \aaa^{\lceil tq \rceil}f_i^q}{K_i}
$$
It is clear that $\aaa^{\lceil t'q \rceil } \subseteq K_i :f_i^q$, so that the length of this term is bounded above by the length of $\aaa^{\lceil tq \rceil}/\aaa^{\lceil t'q \rceil }$. Thus,
$$
e_{HK}(\aaa^{t'}, I) - e_{HK}(\aaa^t, I) \le n \lim_{q\rightarrow 0} \frac{\lambda(R/\aaa^{\lceil t'q \rceil }) -\lambda(R/\aaa^{\lceil tq \rceil})}{q^d}
$$
$$
=ne(\aaa)(t'^d -t^d)
$$
where the last equality is from Lemma ~\ref{realmult}.
\end{proof}

We now show how the joint Hilbert-Kunz multiplicity is related to tight closure, integral closure, and $\aaa$-tight closure. The result pertaining to $\aaa$-tight closure, Proposition ~\ref{hkchar} is an analog of testing tight closure via Hilbert-Kunz multiplicities (cf. \cite{HH1}, Theorem 8.17), and testing for integral closure via Hilbert-Samuel multiplicities (cf. \cite{NR}).
\begin{prop}
If $I\subseteq J\subseteq I^*$ and $\aaa \subseteq \bbb \subseteq \overline{\aaa}$, then
$e_{HK}(\aaa^t, I)=e_{HK}(\bbb^t, J)$ for all $t>0$.
\end{prop}
\begin{proof}
The hypothesis implies that there exists $h \in R^o$, such that
$h \bbb^{\lceil tq\rceil}  J\q \subseteq \aaa^{\lceil tq\rceil} I\q$ (if $\bbb = (b_1, \ldots, b_n)$, for each $b_i$ there exists $h_i \in R^o$ such that $h_i b_i^n \in \aaa ^n $ for all $n$; $\bbb^{\lceil tq\rceil}$ is generated by $b_1^{i_1} \cdots b_n^{i_n}$ with $i_1 + \ldots + i_n = \lceil tq \rceil $, and choosing $h'= h_1 \cdots h_n$ we have $h' b_1^{i_1} \cdots b_n^{i_n} \in \aaa^{i_1} \cdots \aaa^{i_n} = \aaa^{\lceil tq \rceil }$.)
It follows that
$$
\lambda\left( \frac{\bbb ^{\lceil tq \rceil}J\q}{\aaa ^{\lceil tq\rceil} I\q}\right)\le
\lambda \left( \frac{\aaa ^{\lceil tq\rceil} I\q :h}{\aaa ^{\lceil tq\rceil} I\q}\right) = \lambda \left( \frac{R}{(\aaa ^{\lceil tq\rceil} I\q, h)}\right),
$$ which is a joint Hilbert-Kunz function over $R/h$, and thus bounded by $\mathcal{O}(q^{d-1})$.
\end{proof}
\begin{prop}\label{hkchar}
Assume that $R$ has test elements for the usual tight closure. Let $\aaa, I, J$ be $\m$-primary ideals in $R$ and let $t>0$. Assume $I\subseteq J$.
Then $J\subseteq {}^{\aaa^t}I^*$ if and only if $e_{HK}(\aaa^t, I)=e_{HK}(\aaa^t, J)$.
\end{prop}
\begin{proof}
Assume that $J\subseteq {}^{\aaa^t}I^*$, so that $c\aaa^{\lceil tq\rceil} J^q\subseteq \aaa^{\lceil tq\rceil} I\q$ for some $c \in R^o$. Then 
$$
\lambda\left(\frac{\aaa^{\lceil tq \rceil} J\q}{\aaa^{\lceil tq \rceil} I\q}\right) \le \lambda\left(
\frac{\aaa^{\lceil tq \rceil} I\q :c}{\aaa^{\lceil tq \rceil} I\q}\right)=
\lambda\left( \frac{R}{(c, \aaa^{\lceil tq \rceil} I\q)}\right),
$$
which is a mixed Hilbert-Kunz function over the $d-1$ dimensional ring $R/c$, and therefore it is bounded by $\mathcal{O}(q^{d-1})$. This shows that $e_{HK}(\aaa^t, I)=e_{HK}(\aaa^t, J)$.

Conversely, assume that $e_{HK}(\aaa^t, I)=e_{HK}(\aaa^t, J)$.
Fix an element $x \in J$.  We want to show $x \in {}^{\aaa^t}I^*$.

Fix $q_0=p^{e_0}$, and fix generators $g_1, \ldots, g_m$ for $\aaa ^{\lceil tq_0\rceil }$. Let $f$ denote the product of a minimal set of generators for $\aaa$, chosen in $R^o$.
Note that $\lceil tq_0q \rceil \le \lceil tq_0 \rceil q$, and thus we have $\aaa^{\lceil tq_0\rceil q} \subseteq \aaa ^{\lceil tq_0 q\rceil} \subseteq ( \aaa ^{\lceil tq_0\rceil })\q:f^q$ (the last inclusion is Observation  ~\ref{computation} ({\bf a.})).

For each $g_i$, we have
$$
\lambda \left( \frac{(\aaa ^{\lceil tq_0q\rceil } I^{[q_0q]}, \aaa ^{\lceil tq_0q\rceil } x^{q_0q})}{\aaa ^{\lceil q_0q \rceil } I^{[q_0q]}}\right) \ge \lambda \left( \frac{(\aaa ^{\lceil tq_0q\rceil } I^{[q_0q]}, g_i ^{q} x^{q_0q})}{\aaa ^{\lceil tq_0q\rceil } I^{[q_0q]}}\right)
=\lambda \left( \frac{R}{\aaa ^{\lceil tq_0q\rceil } I^{[q_0q]}: g_i^qx^{q_0q}}\right) 
$$
$$
\ge \lambda \left( \frac{R}{{\aaa ^{\lceil tq_0\rceil }}^{[q]}I^{[q_0q]}:(fg_ix^{q_0})^q}\right)=\lambda \left(\frac{(\aaa ^{\lceil tq_0\rceil }I^{[q_0]}, fg_ix^{q_0})\q }{(\aaa ^{\lceil tq_0\rceil }I^{[q_0]})\q }\right).
$$
On the other hand, our assumption implies that 
$$\lambda \left( \frac{(\aaa ^{\lceil tq_0q\rceil } I^{[q_0q]}, \aaa ^{\lceil tq_0q\rceil } x^{q_0q})}{\aaa ^{\lceil tq_0q\rceil } I^{[q_0q]}}\right)
\le \lambda\left( \frac{\aaa^{\lceil tqq_0\rceil }J^{[qq_0]}}{\aaa^{\lceil tqq_0\rceil }I^{[qq_0]}}\right)$$ is bounded above by $\mathcal{O}(q^{d-1})$.
Since $q_0$ is fixed, Theorem 8.17 in \cite{HH1} implies $f g_i x^{q_0} \subseteq (\aaa ^{\lceil tq_0\rceil }I^{[q_0]})^*$.

 Since $g_i$ ranges through the generators of $\aaa ^{\lceil tq_0\rceil }$, we have
$ f\aaa ^{\lceil tq_0\rceil } x^{q_0} \subseteq (\aaa ^{\lceil tq_0\rceil }I^{[q_0]})^*$.
 But this is true for all $q_0$; if we let $c\in R^o$ be a test element for $R$, we get
$$
cf \aaa ^q x^q \subseteq \aaa ^{\lceil tq \rceil }I\q
$$
for all $q>>0$. Since $f\in R^o$, this gives the desired conclusion.
\end{proof}

 We end this section with some formulas relating the joint Hilbert-Kunz multiplicity to the usual Hilbert-Kunz multiplicity and Hilbert-Samuel multiplicity.

\begin{theorem}\label{bounds}
Assume that $\mathrm{dim}(R)\ge 1$.
Let $\aaa, I \subset R$ be $\m$-primary ideals.

\rm{a.} For all $t>0$ we have
\begin{equation}
e_{HK}(\aaa^{t}, I)\le e_{HK}(I) + \frac{ \ell e(\aaa)t^d}{d!},
\end{equation}
where $\ell $ denotes the $*$-spread of $I$, i.e the minimal number of generators of an ideal $J$ minimal with respect to the condition $J^* \supseteq I$.

\rm{b.} If we assume that $R$ is excellent and analytically irreducible then there exists a $t_0>0$ such that the inequality in part a. is equality for all $0\le t \le t_0$.

\rm{c.}
$$ \lim_{t \rightarrow \infty} \frac{ e_{HK}(\aaa^t, I)}{t^d}= \frac{e(\aaa)}{d!}.
$$
In particular, if $\ell >1$ then the inequality in part a. is strict for $t \gg 0$.
\end{theorem}
\begin{proof}
First note that we can replace $I$ by any ideal $J$ with $J \subseteq I \subseteq J^*$ without affecting the result. Thus, we may assume that 
$I=(f_1, \ldots, f_{\ell})$, where $f_1, \ldots, f_{\ell}$ are $*$-independent, i.e. none of them is in the tight closure of the ideal generated by the others. We can also choose all $f_i \in R^o$ by prime avoidance.

We have a filtration
$$
\aaa ^{\lceil tq\rceil}  I \q \subseteq (\aaa ^{\lceil tq\rceil}  I \q, f_1^q) \subseteq \ldots \subseteq (\aaa ^{\lceil tq\rceil}  I \q, f_1^q, \ldots, f_{\ell-1}^q) \subseteq I\q,
$$
and therefore we have
\begin{equation}\label{sum}
\lambda\left( \frac{ I\q } {\aaa^{\lceil tq\rceil}  I\q} \right)
= \Sigma _{i=1}^{\ell} \lambda \left( \frac{R}{ (\aaa ^{\lceil tq\rceil}  I \q, f_1^q, \ldots, f_{i-1}^q): f_i^q}\right)
\end{equation}
Since the denominator in each term in the right hand sum contains $\aaa ^{\lceil tq \rceil }$, the inequality follows by Lemma ~\ref{realmult}.

The second statement follows from Theorem 3.5 (a) in \cite{indep}.

In order to see the last statement, it is enough to restrict to integer exponents $t$. Note that the denominators appearing in the terms on the right hand side of Equation ~\ref{sum} contain $(\aaa^{ tq }, f_1^q, \ldots, f_{i-1}^q)$, and thus
$$
\lambda\left( \frac{ I\q } {\aaa^{ tq}  I\q} \right)
\le \Sigma _{i=1}^{\ell} \lambda \left( \frac{R}{ (\aaa ^{ tq}, f_1^q, \ldots, f_{i-1}^q)}\right)
$$
Consider $i >1$. 
We have
$$
\lim_{q\rightarrow \infty} \lambda\left( \frac{R}{ (\aaa ^{ tq}, f_1^q, \ldots, f_{i-1}^q)}\right)/q^d \le e_{HK}((\aaa^t,  f_1, \ldots, f_{i-1}) )
$$
$$\le \lambda\left(\frac{R}{(\aaa^t,  f_1, \ldots, f_{i-1})} \right) e_{HK}(\m)
$$
(the last inequality follows by taking a filtration of $R/(\aaa^t,  f_1, \ldots, f_{i-1})$ with quotients equal to $R/\m$; also see Lemma 4.2 in \cite{WY}).
As a function of $t$, $\displaystyle \lambda\left(\frac{R}{(\aaa^t,  f_1, \ldots, f_{i-1})} \right)$ is a Hilbert-Samuel function over the ring $R/(f_1, \ldots, f_i)$, which has Krull dimension less than $d$, and therefore dividing by $t^d$ and taking the limit when $t \rightarrow \infty$ yields a limit equal to zero for each of the terms corresponding to $i >1$ in Equation ~\ref{sum}.

Thus we have
$$
\limsup_{t\rightarrow \infty} \frac{e_{HK}(\aaa^t, I)-e_{HK}(I)}{t^d}= \limsup_{t\rightarrow \infty} \lim_{q \rightarrow \infty} \lambda\left( \frac{R}{\aaa^{tq}I\q:f_1^q}\right)/t^dq^d 
$$
$$
\le \lim_{t, q\rightarrow \infty}\lambda\left( \frac{R}{\aaa^{tq}}\right)/t^dq^d =\frac{e(\aaa)}{d!}.
$$
On the other hand, we have $\aaa^{tq} I\q \subset \aaa^{tq}$, and thus
$e_{HK}(\aaa^t, I)\ge t^de(\aaa)/d!$. This proves the equality in part c.

\end{proof}

The following provides a concrete example where part b. in Theorem ~\ref{bounds} works with $t_0=1$.
\begin{example}
Assume $(R, \m)$ is a Cohen-Macaulay ring, and let $\aaa = I =(x_1, \ldots, x_d)$ be generated by a regular sequence.
If $0 \le t \le 1$, then
$$e_{HK}(\aaa^t , I)=\frac{t^d e(\aaa)}{(d-1)!} + e(\aaa)
$$

\end{example}
\begin{proof}
For this choice of $\aaa $ and $I$, each term in the sum on the right hand side of Equation ~\ref{sum} for $0 \le t \le 1$ is equal to $\lambda (R/\aaa^{\lceil tq\rceil} )$, and therefore
$$
\lambda \left( \frac{R}{\aaa^{\lceil tq \rceil} I\q}\right) = d \lambda \left( \frac{R}{\aaa^{\lceil tq\rceil} }\right) + \lambda \left( \frac{R}{I\q}\right).
$$ 
\end{proof}
\section{Test ideals}
The main result of this section, Theorem ~\ref{testideals} shows that the test ideal for the new version of $\aaa$-tight closure coincides with the test ideal for the Hara-Yoshida $\aaa$-tight closure for $R_{_+}$-primary ideals $\aaa$ in a graded Gorenstein ring.
\begin{lemma}\label{irred}
Let $(R, \m)$ be a local approximately Gorenstein ring of characteristic $p>0$. Let $\{I_t\}$ be a sequence of $\m$-primary irreducible ideals, such that for every $k$ there exists $t$ with $\m ^k \subseteq I_t$.

Then
$\tau(\aaa)=\cap_t (I_t : I_t^{* \aaa})$
and $ T_{\aaa}=\cap_t (I_t : {}^{\aaa}I_t ^*)$.
\end{lemma}
\begin{proof}
We'll prove the second statement (the proof for the first one is slightly easier). The inclusion $T_{\aaa}\subseteq \cap (I_t : {}^{\aaa}I_t ^*)$ is clear by definition.
Consider $c \in \cap (I_t : {}^{\aaa}I_t ^*)$.
First we show that $c ({}^{\aaa}I ^*)\subseteq I$, where $I$ is an arbitrary $\m$-primary ideal. The assumption guarantees that there exists $t$ such that $I_t \subseteq I$, and since $I_t$ is irreducible, we can write $I= I_t :K$ for some ideal $K$. Let $x \in {}^{\aaa}I ^*$.
Then there exists $d \in R^o$ such that $dx^q \aaa^q \subseteq \aaa^q I_t\q \subseteq \aaa^q (I_t\q :K\q) \subseteq (\aaa^q I_t \q): K\q$. Thus, $dx^qK\q \aaa^q \subseteq \aaa^q I\q$, which shows that $xK \subseteq {}^{\aaa}I_t ^*$. We have $cxK \subseteq I_t$ by the choice of $c$, and thus $cx \in I_t : K =I$.

Now consider $I$ an arbitrary ideal. We can write $I = \cap _n (I + \m^n)$, an intersection of $\m$-primary ideals. Let $x \in {}^{\aaa}I ^*$. We need to show that $cx \in I$. Note that $x \in {}^{\aaa}(I+ \m^n) ^*$ for all $n$, and therefore $cx \in I+ \m^n$ for all $n$ since we have already proved this for $\m$-primary ideals. Intersecting over all $n$ yields the desired conclusion.
\end{proof}
Throughout the rest of this section, $R$ will be assumed to be a Gorenstein positively graded algebra over a field of Krull dimension $d$ and $a$-invariant $a$. We let $x_1, \ldots, x_d$ be a system of parameters with $\mathrm{deg}(x_i)=\alpha $ for all $i$, and $I_t:=(x_1^t, \ldots, x_d^t)$. Let $u$ denote a homogeneous socle generator for $(x_1, \ldots, x_d)$, i.e. $u \in (x_1, \ldots, x_d):R_{_+}\, \backslash (x_1, \ldots, x_d)$, and let $\delta:=\mathrm{deg}(u)$. Note that $\delta = \alpha d + a$, since the $a$-invariant may be defined as the degree of
$$
\left[ \frac{u}{x_1\ldots x_d}\right] \in H_{R_{_+}}^d(R).
$$
Note that $I_t\q = I_ {tq}$, and its socle is generated by $(x_1\cdots x_d)^{tq-1}u$. We will use $\delta _t$ to denote the degree of the socle generator for $I_t$. More precisely,
$\delta _t =\mathrm{deg}((x_1 \cdots x_d)^{t-1} u) =(t-1)\alpha d + \delta$.
Note that we have $\delta_{ tq }=q\delta_t -(q-1)a$.

Fix the notation established before Proposition ~\ref{highdeg}.
\begin{lemma}\label{degree}
With notations as above, we have
$$
I_t : R_{_+}^N\subseteq I_t + R_{\ge \delta _t -(N-1) \beta}.
$$
Moreover, if $R$ is standard graded, then we have equality.
\end{lemma}
\begin{proof}
We prove the claim by induction on $N$.
 For the case $N=1$, $I_t : R_{_+} =(I_t, (x_1\cdots x_d)^{t-1}u)\subseteq I_t + R_{\ge \delta _t}$ by the definition of $\delta _t$.

To see that the other inclusion holds in the standard graded case, note that every homogeneous element not in $I_t$ must have a multiple in the socle of $I_t$, and thus must have degree $\le \delta _t$.

Assume the claim is true for $N-1$. Note that 
$I_t : R_{_+}^N=(I_t :R_{_+}^{N-1}): R_{_+}$. By the induction hypothesis, we can write $I_t :R_{_+}^{N-1}=(I_t, v_1, \ldots, v_r)$ with $\mathrm{deg}(v_i)\ge \delta _t-(N-2)\beta$ for all $i$. If $v \in (I_t :R_{_+}^N)\, \backslash \, I_t :R_{_+}$, then we have $vy_j\equiv a_1v_1+ \ldots + a_rv_r$ (mod $I_t$) for some $1 \le j \le s$, where $a_i\in R$ can be assumed homogeneous and not all zero. Thus, $\mathrm{deg}(v) + \mathrm{deg}(y_j) \ge \mathrm{deg}(v_i)$ for some $i$, and the desired inclusion follows.

For the other inclusion in the standard graded case: if $x \in R_{\ge \delta _t -( N-1)}$, then for all $y \in R_{_+}$ we have $xy \in R_{\ge \delta_t -(N-2)}\subseteq I_t :(R_{_+})^{N-1}$ by the induction hypothesis.
\end{proof}
\begin{theorem}\label{testideals}
Let $R$ be a Gorenstein finitely generated graded algebra over a field of positive characteristic. Assume that the Krull dimension $d$ of $R$ is at least 2. Let $\aaa$ be a homogeneous ideal which is primary to $R_{_+}$. Then $\tau (\aaa)=T_{\aaa}$.
\end{theorem}
\begin{note}
The statement of the theorem is not true if the Krull dimension is $d=1$, since then we can take $\aaa =(f)$ to be a principal ideal, and we have ${}^{\aaa}I^*=I^*$ and $I^{*\aaa}=I^*:f$ for every ideal $I$. It follows that $T_{\aaa}=\tau$, and
$\tau({\aaa}) = \bigcap_I(I:(I^*: f))= \bigcap_I(I:(I: \tau f))=\bigcap_I(I, \tau f) =\tau f$, where the intersection is taken over all parameter ideals $I$ (see Lemma~\ref{irred}).
\end{note}

\begin{proof}
Fix $c \in R^o$ a homogeneous element such that for all ideals $I\subset R$ we have
$x \in I^{*\aaa} \Rightarrow cx^q\aaa^q \subseteq I\q$. Such a $c$ is called a test element for $\aaa$-tight closure, and the existence of such an element is guaranteed by Theorem 1.7 in \cite{HY}. Fix $k \ge l$ integers such that $R_{_+}^k \subseteq \aaa \subseteq R_{_+}^l$. 

 Due to Lemma ~\ref{irred}, it is enough to prove that $I_t^{*\aaa}={}^{\aaa}I_t^*$ for all $t\gg 0$. Since both $I_t^{*\aaa}$ and ${ }^{\aaa}I_t^*$ are homogeneous ideals, Proposition ~\ref{highdeg} implies that it is enough to show that $x \in I_t^{*\aaa} \Rightarrow \mathrm{deg}(x)\ge k\beta -l\beta' + t\alpha $ when $t \gg 0$. 

Since $R_{_+}^k\subseteq \aaa$, we have $x \in I_t^{*\aaa} \Rightarrow cR_{_+}^{kq}x^q \subseteq I_t\q$. Thus, it follows that $cx^q \in I_{tq} : R_{_+}^{kq}$.

Applying Lemma~\ref{degree}, we see that
 $x \in I_t^{*\aaa}$ implies that for all $q=p^e$ we have either
$cx^q \in I_t\q$, or else
$\mathrm{deg}(c) + q\mathrm{deg}(x)\ge d(tq -1) \alpha + \delta - (kq -1)\beta$.
If $cx^q \in I_t \q$ for all $q=p^e\gg 0$, then $x \in I^* \subseteq {}^{\aaa}I^*$.
Otherwise, it follows that $\mathrm{deg}(x)\ge  dt\alpha  -k \beta$. Since $d>1$, when we choose $t\gg 0$ we have $dt\alpha - k\beta \ge t\alpha + k\beta -l\beta'$, and thus Proposition ~\ref{highdeg} can be applied to show that $x \in {}^{\aaa} I_t^*$.
\end{proof}
We end this section with explicit an computation of test ideals for $\aaa = R_{_+} ^r$, when $R$ is a Gorenstein graded ring. We note that our result is similar to that in Proposition 5.8 in \cite{HY}, but under different assumptions.

\begin{prop}
Let $R$ be a standard graded Gorenstein algebra over a field. With notations as above, we have
$$
 I_t^{*R_{_+}^r} = I_t^* + I_t : R_{_+}^{a+1+\lfloor r \rfloor }
$$
for all $r\ge 0$. Thus, 
$\tau (R_{_+}^r)=\tau (R) \cap R_{_+}^{a+1+\lfloor r \rfloor }$.
\end{prop}

In particular, if $R$ is F-rational, we have $\tau (R)=R$, and thus $\tau (R_{_+}^r)=R_{_+}^{a+1+\lfloor r \rfloor }$.
This also follows from Proposition 5.8 in \cite{HY}, where $R$ is not necessarily graded (instead, F-rationality of the Rees ring $R[R_{_+}t]$ is required).

 Also, the results of \cite{HS} and \cite{Hara} imply that when $R$ is obtained from a characteristic zero ring by reduction to positive characteristic $p \gg 0$, we have $\tau (R)=R_{_+}^{a+1}$, and thus  $\tau (R_{_+}^r)=R_{_+}^{a+1+\lfloor r \rfloor }$ also holds.


\begin{proof}
Let $x \in  I_t^{*R_{_+}^r}$ be a homogeneous element, so that
$cx^q R_{_+}^{\lceil rq \rceil } \subseteq I_t \q=I_{tq}$ for some homogeneous $c \in R^o$. Then
$cx^q \in I_{tq}:R_{_+}^{\lceil rq \rceil }=I_{tq} + R_{\ge \delta _{tq} -\lceil rq \rceil +1}$ be Lemma~\ref{degree}. Thus we either have $cx^q \in I_t\q$ for all $q\gg 0$, in which case $x \in I_t^*$, or else we have
$\mathrm{deg}(c) + q\mathrm{deg}(x) \ge \delta _{tq} -\lceil rq \rceil +1=q\delta _t -(q-1) a -\lceil rq \rceil +1$ for infinitely many values of $q=p^e$.
Dividing each side by $q$ and taking the limits when $q\rightarrow \infty$ yields
$\mathrm{deg}(x) \ge \delta_t - a -r$, and since $\mathrm{deg}(x)$ is an integer, this means
$\mathrm{deg}(x)\ge \delta_t - a -\lfloor r \rfloor$.
For every homogeneous element $y \in R_{a+1 + \lfloor r \rfloor}$, we have $xy \in R_{\ge \delta _t +1}\subset I_t$. This proves $I_t^{*R_{_+}^r}\subseteq I_t^* + I_t : R_{_+}^{a+1+\lfloor r \rfloor }$.

Conversely, consider $x \in I_t : R_{_+}^{a+1+\lfloor r \rfloor }=I_t +R_{\ge \delta _t -a - \lfloor r \rfloor }$. If $x \in I_t$, there is nothing to prove. Otherwise, we have $\mathrm{deg}(x^q)\ge q\delta _t-aq -\lfloor r \rfloor q\ge  q\delta _t-aq-\lceil rq \rceil$. Choosing $c \in R_{\ge a+1}$ yields 
$\mathrm{deg}(cx^q) \ge \delta_{tq}-\lceil rq \rceil +1$, so that
$cx^q R_{_+}^{\lceil rq \rceil }\subseteq  R_{\ge \delta _{tq} +1}\subset I_t \q$, and thus $x \in I_t^{*R_{_+}^r}$. We note that this inclusion can also be obtained as a Corollary of Theorem 2.7 in \cite{HY}.
\end{proof}
\section{Jumping numbers}
The results of this section address the following question:
\begin{question}
Given ideals $\aaa, I\subset R$, and a fixed $t_0 \ge 0$, does there exist an $\epsilon >0$ such that 
$I^{*\aaa^t}=I^{*\aaa^{t_0}}$, and ${}^{\aaa^t}I^*={}^{\aaa^{t_0}}I^*$ for all $t \in [t_0, t_0 + \epsilon]$?
\end{question}
This question is somewhat related to the notion of jumping numbers for test ideals. The jumping numbers are defined to be the positive real numbers $c$ such that $\tau (\aaa ^c) \ne \tau (\aaa ^{c - \epsilon})$ for any $\epsilon >0$ (a similar notion for multiplier ideals has been introduced in \cite{ELSV}). These have been studied extensively in recent research (\cite{MTW}, \cite{BMS}). In our context, if for a given $t_0$ an $\epsilon $ can be found that does not depend on the ideal $I$, it follows that there are no jumping numbers between $t_0$ and $t_0 + \epsilon$.
We give positive answers to our question in several particular cases. A positive answer implies that for a given $I$, $I^{*\aaa^t}$ and ${}^{\aaa^t}I^*$ are constant on intervals of the form $[t_0, t_1)$. We will call $t_1$ a jumping number for the ideal $I$ if $I^{*\aaa^{t_0}}=I^{*\aaa^t}$ for all $t \in[t_0, t_1)$ for some $t_0 < t_1$, but $I^{*\aaa^{t_0}}=I^{*\aaa^{t_1}}$ (or  ${}^{\aaa^t_0}I^*= {}^{\aaa^t}I^*$ but  ${}^{\aaa^{t_0}}I^*\ne {}^{\aaa^{t_1}}I^*$).

The following observation shows that it will be enough to check only one inclusion in order to answer the above question in the affirmative.
\begin{obs}\label{smaller}
Let $I, \aaa \subset R$ be fixed ideals, and $0 \le t < t'$ real numbers. Then $I^{*\aaa^t}\subseteq I^{*\aaa^{t'}}$, and 
${}^{\aaa^t}I^* \subseteq {}^{\aaa^{t'}}I^*$.

However, note that it is not always true that $\aaa \subseteq \bbb \Rightarrow {}^{\bbb}I^*\subseteq {}^{\aaa}I^*$, while the corresponding statement is trivially true for the Hara-Yoshida version.
\end{obs}
\begin{proof}
The statement for the Hara-Yoshida version is trivial, since 
$\aaa^{\lceil t'q \rceil } \subseteq \aaa^{\lceil tq \rceil }$.

Consider $x \in {}^{\aaa^t}I^*$, so that 
$cx^q \aaa^{\lceil tq \rceil } \subseteq \aaa^{\lceil tq \rceil} I\q$. Multiplying each side by arbitrary elements in $\aaa^{\lceil t'q \rceil - \lceil tq \rceil }$ yields the desired conclusion.

For the claim in the last paragraph, take for example $\aaa =(f)$, with $f \in (\bbb):=(x, y)^2$, $I=(x^2, y^2)$, in the ring $R=k[x, y]$. Then ${}^{\bbb}I^*=(x^2, y^2, xy)$, while ${}^{\aaa}I^*=I$.
\end{proof}
\begin{prop}
Assume that $(R, \m)$ is local, and $\aaa, I$ are $\m$-primary ideal. Then for every $t_0\ge 0$, there exists $\epsilon >0$ such that ${}^{\aaa^t}I^*={}^{\aaa^{t_0}}I^*$ for all $t \in [t_0, t_0 + \epsilon]$.
\end{prop}
\begin{proof}
First note that
for each $x \notin {}^{\aaa^{t_0}}I^*$, there exists $\epsilon >0$ such that $x \notin {}^{\aaa^{t_0+\epsilon}}I^*$. This follows from Theorem~\ref{continuity}, and Proposition~\ref{hkchar}, applied to the ideals $I$ and $J=(I, x)$.

Construct a sequence $t_1 > t_2 > \ldots t_n \ldots  > t_0$ recursively as follows:
Choose $x_1 \notin {}^{\aaa^{t_0}}I^*$, and let $t_1>t_0$ such that $x_1 \notin {}^{\aaa^{t_1}}I^*$ (the existence of such a $t_1$ is guaranteed by the previous claim). If $t_1, \ldots, t_k$ have been constructed, then we either have ${}^{\aaa^{t_0}}I^*={}^{\aaa^{t_k}}I^*$, in which case the proof is complete (take $\epsilon = t_k - t_0$), or else we can choose an $x_{k+1} \in {}^{\aaa^{t_k}}I^* \, \backslash {}^{\aaa^{t_0}}I^*$, and, by the previous claim, there exists $t'_{k+1}>t_0$ such that $x_{k+1} \notin {}^{\aaa^{t_{k+1}}}I^*$ (note that we must have $t_k > t_{k+1}$ by Observation~\ref{smaller}). Thus, we have a chain of ideals ${}^{\aaa^{t_0}}I^* \subseteq \ldots {}^{\aaa^{t_n}}I^* \subseteq \ldots \subseteq {}^{\aaa^{t_2}}I^* \subseteq {}^{\aaa^{t_1}}I^*$. Note that the construction of $t_k$ shows that the inclusions are strict unless the recursive process stops. This contradicts the fact that ${}^{\aaa^{t_0}}I^*$ is $\m$-primary.
\end{proof}
The next result deals with the case when $\aaa =(f)$ is a principal ideal, with $f \in R^o$. Note that in this case we only need to consider the Hara-Yoshida version, since ${}^{\aaa^t}I^*=I^*$ for all $t$. It turns out that a positive answer  to the question considered here is related to the existence of test exponents. We review the definition.
\begin{definition}
Let $I\subset R$ be an ideal, and $c \in R^o$ a test element for the usual tight closure. We say that $q_0=p^{e_0}$ is a test exponent for the pair $(I, c)$ if $cx^q \in I\q$ for any one choice of $q \ge q_0$ implies $x \in I^*$.
\end{definition}
Test exponents were introduced in \cite{HHexp}, where it is shown that their existence is closely related to the localization problem for tight closure.

\begin{lemma}\label{whatitmeans}
Assume $\aaa =(f)$ is a principal ideal, with $f \in R^o$. Assume that $R$ has test elements for the usual tight closure. Then 
$x \in I^{*\aaa^t} \Leftrightarrow x^q f^{\lceil tq \rceil } \in (I^{[q]})^*$ for all $q$.

In particular, if there exists $q_1$ such that $tq_1 \in {\bf Z}$, then $x \in I^{*\aaa^t} \Leftrightarrow x^{q_1}f^{tq_1}\in (I^{[q_1]})^*$.
\end{lemma}
\begin{proof}
Note that the following inequalities hold for all $q$:
$$
\lceil tq \rceil -1 < tq \le \lceil tq \rceil < tq +1.
$$
It follows that
$$
q_1 \lceil tq_2 \rceil -q_1 < \lceil tq_1q_2 \rceil < q_1 \lceil tq_2 \rceil +1.
$$
Assume that $x \in I^{*\aaa^t}$, and let $q=q_1q_2$. Then there exists $c \in R^o$ such that 
$cx^{q_1q_2}f^{\lceil tq_1q_2 \rceil } \in I^{[q_1q_2]}$, which implies
$c f x^{q_1q_2} f^{\lceil tq_2 \rceil  q_1} \in I^{[q_1q_2]}$. Since $cf \in R^o$, this shows that 
$x^{q_2} f^{\lceil tq_2 \rceil}\in (I^{[q_2]})^*$ for all $q_2$.

Conversely, assume that $x^{q_2} f^{\lceil tq_2 \rceil}\in (I^{[q_2]})^*$ for some $q_2$, and let $c \in R^o$ be a test element for the usual tight closure.
Then $cx^{q_1q_2}f^{\lceil tq_2 \rceil q_1} \in I^{[q_1q_2]}$, which implies 
$cx^{q_1q_2} f^{\lceil tq_1q_2 \rceil + q_1} \in I^{[q_1q_2]}$
for all $q_1$. Fix $q_1$ and allow $q_2$ to vary 
$cf^{q_1} x^q f^{\lceil tq \rceil } \in I^{[q]}$ for all $q \gg 0$. Since $cf^{q_1}\in R^o$, this shows that $x\in I^{*\aaa^t}$.
\end{proof}

\begin{prop}\label{testexp}
Let $\aaa =(f)$ with $f\in R^o$, $I \subset R$ an arbitrary ideal, and $t_0\ge 0$. Assume that there exists $q_1$ such that $t_0q_1$ is an integer, and assume that $c \in R^o$ is a test element for the closure such that there exists a test exponent $q_0$ for the ideal $I^{[q_1]}$ and the test element $cf^{q_1}$.

Then we have $I^{*\aaa^t}=I^{*\aaa^{t_0}}$ for all $t \in [t_0, t_0 + 1/q_0]$.
\end{prop}
\begin{proof}
Let $\displaystyle t=t_0 + \frac{1}{q_0}$. Let $q=q_1q_0$ so that $tq=t_0q + q_1$ is an integer. Assume that $x \in I^{* \aaa^{t}}$; by Lemma~\ref{whatitmeans}, this implies that $x^{q}f^{t_0q}f^{q_1} \in I^{[q]})^*$. Since $c$ is a test element, we have
$cf^{q_1} (x^{q_1}f^{t_0q_1})^{q_0} \in (I^{[q_1]})^{[q_0]}$. Since $q_0$ is a test exponent, this implies $x^{q_1}f^{t_0q_1}\in
(I^{[q_1]})^*$, and thus $x \in I^{*\aaa^{t_0}}$ by Lemma~\ref{whatitmeans}.
\end{proof}
\begin{corollar}
If $(R, \m)$ is a regular local ring, $\aaa=(f)$ is a principal ideal, and $t_0\ge 0$ is such that $t_0q_1 \in {\bf Z}$ for some $q_1=p^{e_1}$,and $q_0$ is such that $f^{q_1} \notin \m^{[q_0]}$, then there are no jumping numbers for the test ideals $\tau (\aaa^t)$ in the interval $[t_0, t_0 + 1/q_0]$.
\end{corollar}
\begin{proof}
Since $I^*=I$ for every ideal $I$, we can take $c=1$, and note that a $q_0$ with the property that $f^{q_1} \notin \m^{[q_0]}$ is a test exponent for $(I^{[q_1]}, f^{q_1})$ for any ideal $I$. Indeed, if $x\notin I^{[q_1]}$ and $f^{q_1}x^{q_0} \in I^{[q_1q_0]}\Rightarrow f^{q_1} \in I^{[q_1q_0]}:x^{q_0} =(I^{[q_1]}:x)^{[q_0]}\subseteq \m^{[q_0]}$, contradicting the choice of $q_0$.
\end{proof}
Note that if $t_0=0$, then the converse of Proposition ~\ref{testexp} holds, in the sense that a positive answer to the question discussed here implies existence of test exponents for the usual tight closure. Recall that $I^{* \aaa^{0}}=I^{*R}=I^*$ is the usual tight closure.
\begin{prop}
Let $\aaa =(f)$ with $f\in R^o$ a test element for tight closure, $I \subset R$ an arbitrary ideal.

 Assume that $q_0=p^{e_0}$ is such that $I^{*\aaa^{1/q_0}}=I^*$. Then $q_0$ is a test exponent for the pair $(I, f)$.
\end{prop}
\begin{proof}
Let $x \in R$ be such that $fx^{q_1} \in I^{[q_1]}$ for some $q_1\ge q_0$. Then for all $q \ge q_1$ we have $f^{q/q_1}x^q \in I\q$, and therefore $f^{q/q_0}x^q \in I\q$, which shows that $x\in I^{*\aaa^{1/q_0}}=I^*$ by assumption.
\end{proof}

In the next result, $I$ and $\aaa$ are arbitrary ideals, but we restrict attention to $t_0=0$. 
\begin{prop}Assume that $(R, \m)$ is an excellent analytically irreducible local domain.
Let $I, \aaa \subset R$ be ideals. Let $x \notin I^*$. Then there exists $q_0=p^{e_0}$ such that $x \notin I^{*\aaa ^{1/q_0}}$.
\end{prop}
Note that our result is not quite sufficient to give an affirmative answer to the question raised in the beginning of the section for this case, since $q_0$ is allowed to depend on $x$.
\begin{proof}
By Proposition 2.4 in \cite{Ab}, there exists $q_1$ such that $I\q : x^q \subset \m^{[q/q_1]}$ for all $q \ge q_1$.
Assume by contradiction that $x \in {}^{\aaa^{1/q_0}}I^*$ for every $q_0$. This means that
$cx^q \aaa ^{q/q_0} \subset I^{[q]}$ for all $q \gg0$, so that
$c \aaa^{q/q_0}\subseteq I\q : x^q \subseteq \m^{[q/q_1]}$. Let
$q=q_0q_1 Q$. Then we obtain
$c \aaa ^{q_1 Q}\subseteq \m^{[q_0Q]}$, which implies 
$\aaa ^{q_1}\subseteq (\m^{[q_0]})^*$. Since $q_1$ is fixed, this is false for $q_0 \gg 0$.
\end{proof}

At the other end of the spectrum, we ask the following question.
\begin{question}
If $I, \aaa$ are fixed ideals, and $N$ is such that 
${}^{\aaa^N}I^*={}^{\aaa^{N'}}I^*$ for all $N' \ge N$ (note that such an $N$ exists by the Noetherian property) describe ${}^{\aaa^N}I^*$.
\end{question}
We will use ${ }^{\aaa^{\infty }}I^*$ to denote ${}^{\aaa^N}I^*$ when $N$ is as above. Note that a similar definition for the Hara-Yoshida version of $\aaa$-tight closure would yield the whole ring when $I$ is an $\m$-primary ideal, since for $N\gg 0$ we have $\aaa ^N \subseteq I$, and $\aaa ^{\lceil Nq/k\rceil }\subseteq I\q$, where $k$ is the number of generators of $\aaa$.
When $\aaa =(f)$ is a principal ideal with $f$ a non-zerodivisor, we have ${ }^{\aaa^{\infty}}I^*=I^*$ for every ideal $I$.

We always have ${}^{\aaa ^{\infty}} I^* \subseteq \overline{I}$. Observation ~\ref{isic} implies that when $R$ is standard graded, $\aaa = R_{_+}^r$ for some $r>0$, and $I$ is homogeneous with all generators of the same degree, we have
${}^{*\aaa^{\infty}}I^*=\overline{I}$. However, Example~\ref{notic} shows that for $R=k[x, y]$, $I=(x^2, y^4)$, and $\aaa=(x, y)$ we have
${}^{\aaa^{\infty} }I^*\ne \overline{I}$. In fact in this example it is easy to check that ${}^{\aaa^{\infty} }I^*=(x^2, y^4, xy^3)$. More generally, we note the following:
\begin{prop}
Let $R$ be a two-dimensional standard graded normal domain,
 let $I=(f_1, f_2)$ be a homogeneous parameter ideal. Let $d=\mathrm{max}(\mathrm{deg}(f_1), \mathrm{deg}(f_2))$.
Then
$ {}^{R_{_+}^{\infty}}I^*=I +R_{\ge d}$.
\end{prop}
\begin{proof}
Say that $d =\mathrm{deg}(f_1)$.

Let $x \in {}^{R_{_+}^{\infty}}I^*$, so $x \in {}^{R_{_+}^n}I^*$ for some $n$. Assume that $\mathrm{deg}(x) <d$. For some homogeneous $c \in R^o$, and for all $y$ among a minimal set of generators of $R_{_+}^{nq}$ we have
$cx^q y =a_1f_1^q + a_2f_2^q$ with $a_1, a_2 \in  R_{_+}^{nq}$. If $\mathrm{deg}(x) <d$, it follows by comparing degrees that$cx^q y = a_2f_2^q$, so that $c R_{_+}^{nq}x^q \subseteq R_{_+}^{nq}f_2^2$. But this implies that $x \in \overline{(f_2)}=(f_2)$.

Conversely, assume that $\mathrm{deg}(x) \ge d$. For $n \gg 0$ we have $x\in I^{*R_{_+}^n}$, so that there exists $c\in R^o$ such that for all $y \in R_{_+}^{nq}$,
$cx^q y =a_1f_1^q + a_2f_2^q$ with $a_1, a_2 \in R$. Comparing degrees, we see that $a_1, a_2 \in R_{_+}^{nq}$.
\end{proof}


\begin{thebibliography}{99}

\bibitem[Ab]{Ab} I. Aberbach, {\em Extensions of weakly and strongly F-regular rings by flat maps}, J. Algebra, {\bf 241} (2001), 799-807.

\bibitem[Bt]{Bt} P. B. Bhattacharya, {\em  The Hilbert function of two ideals},  Proc. Cambridge Philos. Soc.  {\bf 53}  (1957), 568--575.

\bibitem[BMS]{BMS} M. Blickle, M. Mustata, and K. E. Smith, {\em Discretness and rationality of F-thresholds}, preprint.

\bibitem[ELSV]{ELSV} L. Ein, R. Lazarsfeld, K. E. Smith and D. Varolin, {\em Jumping coefficients for multiplier ideals}, Duke Math. J., {\bf 123} (2004), 469--506.

\bibitem[Ep]{Ep} N. Epstein, {\em A tight closure analogue of analytic spread},  Math. Proc. Cambridge Philos. Soc., {\bf 139}(2005), no.2, 371--383.



\bibitem[Hn]{Hn} D. Hanes, {\em Notes on the Hilbert-Kunz function},  J. Algebra {\bf 265} (2003), no. 2, 619--630.


\bibitem[H1]{Hara} N. Hara, {\em Characterization of rational singularities in terms of the injectivity of Frobenius}, Amer. J. Math.  {\bf 120}  (1998),  no. 5, 981--996. 

\bibitem[H2]{Ha} N. Hara, {\em Geometric interpretation of test ideals},  Trans. Amer. Math. Soc.  {\bf 353}  (2001),  no. 5, 1885--1906.

\bibitem[HT]{HT} N. Hara and S. Tagaki, {\em On a generalization of test ideals},
Nagoya Math. J. {\bf 175} (2004), 59--74.
 
\bibitem[HY]{HY} N. Hara and K. Yoshida, {\em A generalization of tight closure and multiplier ideals},  Trans. Amer. Math. Soc.  {\bf 355}  (2003),  no. 8, 3143--3174.


\bibitem[HH1]{HH1} M. Hochster, C. Huneke,  {\em Tight closure, invariant theory, and the Brian\c con-Skoda theorem}, J. Amer. Math. Soc. {\bf 3} (1990), 31--116.


\bibitem[HH2]{HHexp} M. Hochster and C. Huneke, {\em Localization and test exponents for tight closure},  Michigan Math. J.  {\bf 48}  (2000), 305--329. 
\bibitem[HS]{HS} C. Huneke and K. Smith, {\em Tight closure and the Kodaira vanishing theorem},
J. Reine Angew. Math. {\bf 484} (1997), 127--152. 


\bibitem[Mo]{Mo} P. Monsky, {\em The Hilbert-Kunz function}, Math. Ann., {\bf 263} (1983), 43--49.

\bibitem[MTW]{MTW} M. Mustata, S. Tagaki, and K.-i. Watanabe, {\em F-thresholds and Bernstein-Sato polynomials}, European Congress of Mathematics, 341--364, Eur. Math. Soc., Zrich, 2005.


\bibitem[NR]{NR} D. G. Northcott, and D. Rees,
{Reductions of ideals in local rings},
Proc. Cambridge Philos. Soc. {\bf 50}, (1954), 145--158. 

\bibitem[S1]{S1} K. E. Smith, {\em Tight closure in graded rings}, J. Math. Kyoto Univ., {\bf 37} no. 1 (1997), 35--53.
 
\bibitem[S2]{Sm} K. E. Smith, {\em The multiplier ideal is a universal test ideal, Comm. Algebra} {\bf 28} (2000), 5915--5929.


\bibitem[Vr]{indep} A. Vraciu, {\em *-Independence and special tight closure}, J. Algebra, {\bf 249} (2002), 544--565.

\bibitem[WY]{WY} K. Watanabe and K. Yoshida, {\em Hilbert-Kunz multiplicity and an inequality between multiplicity and colength}, J. Algebra, {\bf 230} (2000), 295--317. 
\end{thebibliography}
\end{document}